\documentclass[a4paper]{article}
\usepackage{amsmath,amssymb,amsthm,epsfig}
\usepackage{graphicx}
\usepackage[utf8]{inputenc}
\usepackage{color}
\usepackage[pdftex,colorlinks]{hyperref}
\usepackage{subfigure}

\usepackage[all]{xy}

\makeatother

\newtheorem{theor}{theorem}[section]
\newtheorem{theorem}[theor]{Theorem}
\newtheorem{defi}[theor]{Definition}

\newtheorem{corollary}[theor]{Corollary}

\newtheorem{question}[theor]{Question}

\newcommand{\CD}{\operatorname{CD}}

\newcommand{\CAT}{\operatorname{CAT}}

\newcommand{\Scal}{\operatorname{Scal}}

\newcommand{\dist}{\operatorname{dist}}

\title{Scalar Curvature on some Open $3$-Manifolds}

\author{G.~Besson}

\begin{document}

\maketitle
\abstract{}

\tableofcontents

\section{Introduction}\label{intro:scalar}

The {\sl Notices of the AMS}, in their volume 58, number 9 of October 2011 (p.1226--1249), presented a beautiful collection of essays in the memory of Shiing-Shen Chern. The contributions are of the highest quality and I was struck by the essay written by the late Manfredo do Carmo entitled \lq\lq On collaborating with Chern". Manfredo do Carmo explained how a course given in 1968, in Berkeley,  by S.~S.~Chern has deeply influenced his career. The course was about a (at the time) recent preprint by J.~Simons which eventually was published as the article \cite{Simons} and dealt with minimal varieties in Riemannian manifolds. According to Manfredo do Carmo the course was given in such an exciting way that he became convinced that studying minimal surfaces was his \lq\lq home". I was lucky to meet several times Manfredo do Carmo in Brazil and not only he was a renowned specialist of minimal varieties but he also had built a strong school and there is certainly very few mathematical departments in Brazil  that does not have a researcher who is a descendant of Manfredo. It seems not exaggerated to say that this grew out of this course given by S.~S.~Chern in Berkeley in the winter of 1968 ! Very few mathematicians can have such a worldwide influence and S.~S.~Chern was one of them.

The paper by Jim Simons became immediately famous. In its first part the author proved what is now called the {\sl Simons equation}, which is an elliptic equation satisfied by the second fundamental form of any minimal variety in any Riemannian manifold. The minimal surfaces (or submanifolds) theory has produced lots of beautiful examples and is still very active nowadays. They can be studied for themselves or for their relations with the topology and the geometry of  the ambiant Riemannian manifold. There is a notion of {\sl stability} related  to the {\sl Jacobi operator}, and {\sl grosso modo} saying that under normal local (compact) deformations of the minimal submanifold the variation of the volume has non negative hessian. The interested reader is referred to  \cite{Col-Min} for the details. 

A striking fact is that stable minimal surfaces in $3$-manifolds have a strong impact on their topology and geometry. This was noticed by several authors, in particular R.~Schoen and S.~T.~Yau (\cite{SchoenYau:stable}) and M.~Gromov and B.~Lawson (\cite{GromovLawson}). In the present text we present a survey of some recents results on a family of manifolds: the contractible open $3$-manifolds, the first example of which was constructed by J.~H.~C.~Whitehead. Stable minimal surfaces are used to prove a rigidity result asserting that, among all of them, $\textbf{R}^3$ is distinguishable by the geometry it can carry. This is by no mean exhaustive but concerns striking examples of (open) $3$-manifolds.

\thanks{I would like to thank Xiaonan Ma for giving me the opportunity to advertise for this topic, and for his patience, and Jian Wang for interesting discussions.}

\section{Stable Minimal Surfaces and Positive Scalar Curvature}\label{sec:stable}

Let us start this short discussion by Theorem 10.2 of \cite{GromovLawson} which we recall here

\begin{theorem}[\cite{GromovLawson} p. 172]\label{theo:GLmain}
Let S be a compact stable minimal surface in a $3$-manifold $X$ with scalar curvature satisfying $\Scal \geq \kappa_0 > o$. Let $\Omega\subset\Sigma$ be a compact connected domain and let $\rho > o$ be a number
such that
  \begin{enumerate}
      \item $\Omega_\rho $ does not meet $\partial\Sigma$.
      \item Image$[H_1(\Omega )\rightarrow H_1(\Omega_\rho )]\ne 0$.
  \end{enumerate}
Then,
$$\rho <\frac{\pi}{\sqrt{\kappa_0}}\,.$$
\end{theorem}
Here, we denote by $\Omega_\rho :=\{x\in \Sigma : d_\Sigma (x, \Omega )\leq \rho \}$. The proof of this result given on page 178 of \cite{GromovLawson} uses in an essential way the stability of $\Sigma$. This notion of stability implies that for every smooth function $f$ with compact support in $\Sigma\setminus\partial\Sigma$ one has,
$$\int_\Sigma\Big\{ \vert\nabla f\vert^2+Kf^2-\frac{1}{2}\kappa_0 f^2 \Big\}dA\geq 0\,,$$
where $K$ is the Gaussian curvature of $\Sigma$ and $dA$ the area element of the metric on $\Sigma$ induced from the one on $X$.

A consequence of Theorem \ref{theo:GLmain}, stated as Theorem 10.7, is the following statement.

\begin{theorem}\label{theo:GLcompact}
 Let $X$ be a compact $3$-manifold, possibly with a non empty boundary, and suppose that $X$ is endowed with a metric of scalar curvature greater or equal to $1$. Then any closed curve $\gamma\subset X$ such that
  \begin{enumerate}
    \item $[\gamma]=0$ in $H_1(X, \partial X)$,
    \item $\dist (\gamma , \partial X)>2\pi$
  \end{enumerate}
must already bound in its $2\pi$-neighbourhood.  
\end{theorem}

The proof amounts to applying Theorem \ref{theo:GLmain} to the surface $\Sigma$ of least area spanning $\gamma$ modulo $\partial X$, which may intersect $\partial X$. The argument though is not that straightforward and the interested reader can read it on page 177 of \cite{GromovLawson}. 

The upshot of these nice results is the following corollary (Corollary 10.9 in \cite{GromovLawson}) which is relevant to the purpose of this text. On a non compact manifold $X$, we say that the scalar curvature is \emph{uniformly positive} if there exists $\kappa_0>0$ such that, $\Scal\geq\kappa_0$ everywhere on $X$.

\begin{corollary}\label{cor:GLcompact}
A complete $3$-manifold $X$ of uniformly positive scalar curvature and with finitely generated fundamental group is simply-connected at infinity.
\end{corollary}

This corollary makes sense for open $3$-manifolds and it provides an important topological consequence of uniformly positive scalar curvature. Although $X$, in Corollary \ref{cor:GLcompact}, may be non compact the result  is a consequence of Theorem \ref{theo:GLcompact}. For the sake of completeness we end this section by recalling the definition of \emph{simply-connectedness at infinity}.

\begin{defi}
A connected, locally compact space $M$ is  \emph{simply-connected at infinity} if for any compact set $K\subset M$, there exists a  compact set $K'$ containing $K$ so that the induced map $\pi_{1}(M\setminus K')\rightarrow \pi_{1}(M\setminus K)$ is trivial. 
\end{defi}

In other words, for any compact set $K$, the loops which are far away from $K$, say in the complement of $K'$, can be contracted in the complement of $K$. Notice that this notion is preserved by homeomorphisms.

\section{Contractible Open $3$-Manifolds}

\begin{figure}[!h]
\centering\includegraphics[scale=0.8]{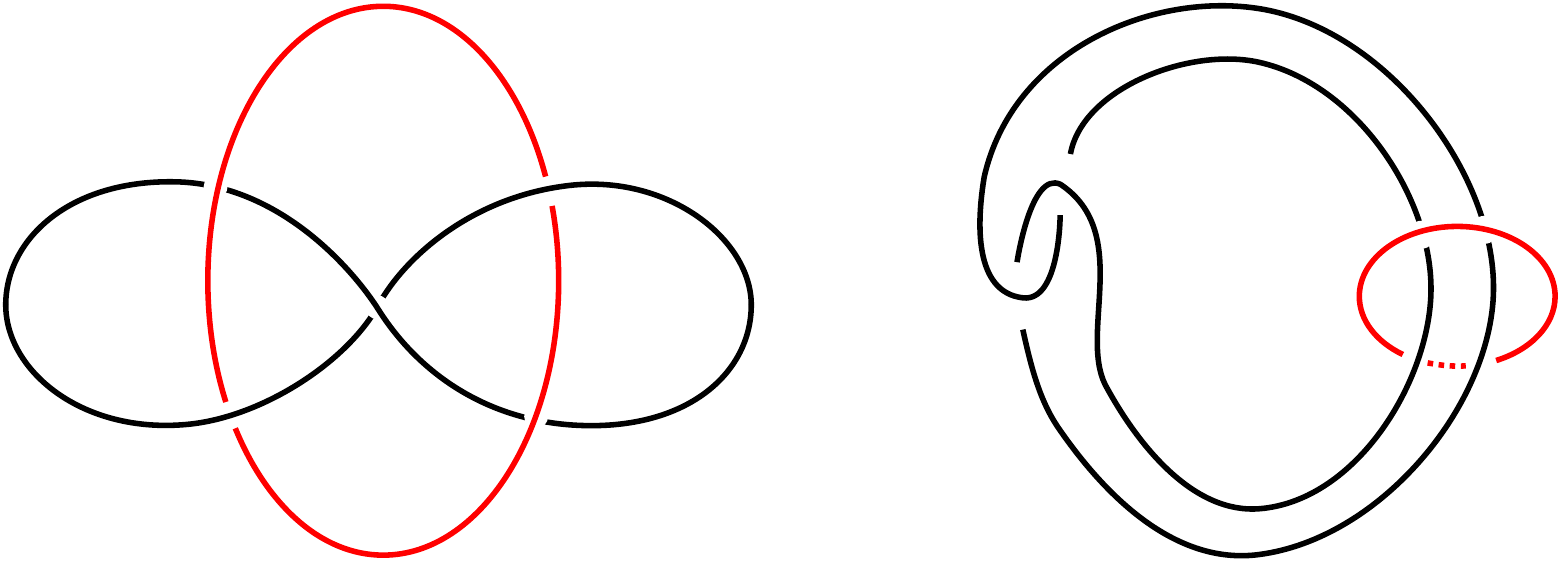}
\caption{Whitehead Link}
\label{Fig1}
\end{figure}

We now describe a series of examples, the contractible open $3$-manifolds which are not homeomorphic to $\textbf R^3$, for which stable minimal disks play an important role. They are trivial in the sense that they are homotopically trivial and their topology is hidden in their structure at infinity. The very first one was introduced by J.~H.~C.~Whitehead in \cite{White1} in an unsuccessful attempt to prove the Poincaré conjecture; let us recall its construction. We start with the Whitehead link which is a link with two components illustrated in Figure \ref{Fig1} below in two different ways. Notice that this link is symmetric; this means that there is an isotopy of the ambient space which reverses the roles played by the black and red curves. 

\begin{figure}[!h]
\centering
\includegraphics[scale=0.9]{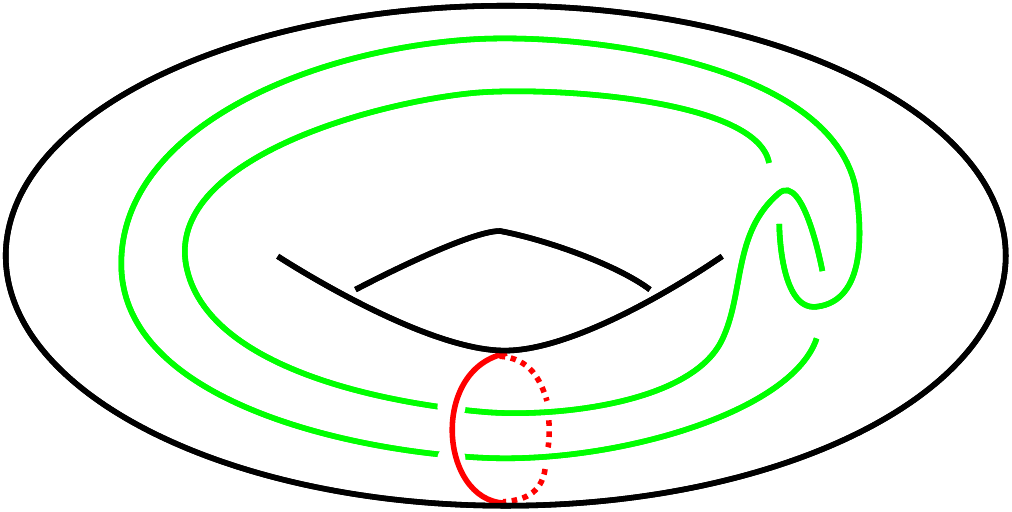}
\caption{$N_1\subset N_2$}
\label{Fig2}
\end{figure}

For a closed solid torus $N$ we define the notion of a meridian curve. A meridian $\gamma\subset \partial N$ is an embedded circle which is nullhomotopic in $N$ but  not contractible in $\partial N$. On Figure \ref{Fig2} the red curve is a meridian.

 We now choose an open solid torus $N_{1}$. It can be viewed as an unknotted solid torus in $\mathbf{S}^{3}$ and it is well known that the complement of $N_1$ in $\mathbb{S}^{3}$ is another solid torus (this time closed). We then embed $N_{1}$ into a second open solid torus $N_{2}$  as a tubular neighbourhood of the green curve shown in Figure \ref{Fig2}. The green and red curves form a Whitehead link. This is the basic pattern of the construction which we will repeat infinitely many times. Precisely, $N_{2}$ is a solid torus which can be embedded, in the same way, into another open solid torus $N_{3}$. We do this infinitely many times and define the {\sl Whitehead manifold}, $Wh:=\cup_{k=1}^{\infty} N_{k}$. It turns out that $Wh$ is an open subset of $\mathbf{S}^{3}$ (hence of $\mathbf{R}^{3}$ too) whose complement is a fractal closed subset of $\mathbf{S}^{3}$ called the {\sl Whitehead Continuum}. The Whitehead manifold is then an open submanifold of $\mathbf{S}^{3}$, {\sl i.e.} a non compact manifold without boundary.

The symmetry of the Whitehead link allows a second construction which describes the Whitehead continuum as the intersection of a descending sequence of closed solid tori. The above one though, with the increasing sequence $N_k$, is much more flexible since it can support variations such as changing the knot at each step $k$, and yields a family of manifolds some of which are even not embedded in $\mathbf{S}^{3}$.
We then say that an open $3$-manifold is genus one if it is the increasing union of open solid tori $N_k$ so that, for each $k$, $\bar N_k\subset N_{k+1}$ (where $\bar N_k$ is the closure of $N_k$) and such that a disc filling a meridian of $N_{k+1}$ intersects the core of $N_k$. In such a generality they were introduced in \cite{McMil} and there are uncountably many of them, some of which are subsets of $\mathbf S^3$, some not (see \cite{K-Mil}). Note that $\mathbf{R}^{3}$ is not genus one but rather genus zero, since it is an increasing union of $3$-balls. The construction can also be made with handlebodies of higher genus \cite{McMil} and the genus can also change at each stage of the construction, the worst case being when this genus goes to infinity. This produces an incredible zoology of contractible pairwise non homeomorphic $3$-manifolds !

Showing that $Wh$ is contractible is easy with the description. Thanks to a theorem due (again) to Whitehead it suffices to show that all homotopy groups are trivial and this follows from the fact that $N_k$ is homotopically trivial in $N_{k+1}$. To prove that it is not homeomorphic to  $\mathbf{R}^{3}$ we simply show that $Wh$ is not simply-connected at infinity. 

We now want to explore the Riemannian geometry of these spaces, the idea being that among all these contractible $3$-manifolds  $\mathbf{R}^{3}$ should be special. The starting point is the article 
\cite{SchoenYau:stable} by R.~Schoen and S.-T.~Yau in which they prove the following theorem.

\begin{theorem}[\cite{SchoenYau:stable}, Theorem 3]
Let $M$ be a complete non compact 3-dimensional manifold with positive Ricci curvature. Then $M$ is diffeomorphic to $\mathbf{R}^{3}$.
\end{theorem}

The key idea relies on showing that there are no stable complete minimal surface in $M$. Indeed, its Jacobi operator is related to the Ricci curvature of the ambient space whose positiveness would give a contradiction. Recently this result was extended by G.~Liu in \cite{Liu} (see Theorem 2) who gets the conclusion that a contractible $3$-manifold cannot carry a complete Riemannian metric with non negative Ricci curvature unless it is $\mathbf{R}^{3}$. The next step brings us to the article by M.~Gromov and B.~Lawson, \cite{GromovLawson}. Particularly to the corollary of the main theorem of their chapter 10, stated above as Corollary \ref{cor:GLcompact}. One consequence is that a contractible $3$-manifold cannot carry a complete metric with uniformly positive scalar curvature unless it is diffeomorphic to $\mathbf{R}^{3}$; indeed, a result by C.~H.~Edwards (\cite{Edw}) combined with the proof of the Poincaré conjecture shows that a contractible open $3$-manifold which is simply-connected at infinity is homeomorphic to $\textbf{R}^3$.  Let us recall that $\textbf{R}^3$ does carry a complete metric of uniformly positive scalar curvature.

Results in the same spirit then appear in \cite{cwy:taming} where it is proved, in particular, that

\begin{theorem}[\cite{cwy:taming}, Theorem 4.4]\label{cwy}
If a non compact contractible 3-manifold M has a complete Riemannian metric with uniformly positive scalar curvature outside a compact set, then it is homeomorphic to $\mathbf{R}^{3}$.
\end{theorem}

Then, a striking result recently announced by Jian Wang gives a definitive answer. 

\begin{theorem}[\cite{Wang}]\label{theo:Wang}
A contractible open $3$-manifold which admits a complete metric of non negative scalar curvature is homeomorphic (hence diffeomorphic) to $\mathbf{R}^{3}$.
\end{theorem}

In \cite{Wang} the above statement is made with the assumption that the scalar curvature is positive. However a nice argument by J.~Kazdan then allows Jian Wang to state the result as above. Theorem \ref{theo:Wang} is a follow-up of two previous versions solving the same question for subfamilies of contractible open $3$-manifolds: see \cite{Wang1} for the genus 1 case and \cite{Wang2} for the case when the fundamental group at infinity is trivial. It is worth recalling that $\textbf{R}^3$ does carry a metric with positive scalar curvature and therefore Theorem \ref{theo:Wang} appears as a beautiful rigidity result.

The method of proof pertains to the same philosophy than in \cite{SchoenYau:stable} and \cite{GromovLawson}. Let us describe some of the ideas contained in the proof of \ref{theo:Wang}, in the case where the scalar curvature is supposed to be positive and the space is a genus one contractible manifold.

 We call $X$ the open contractible $3$-manifold. With the notations used to define genus one $3$-manifolds, let us consider a meridian $\gamma\subset\partial {\bar N}_k$ (where $\bar N_k$ is the closure of $N_k$), this is an embedded curve which can be filled by a minimizing disk $D_k$. Let us assume, for a while, that this disk is included in $\bar N_k$. Jian Wang shows that the number of connected components of $D_k\cap N_2$ intersecting $N_1$ goes to infinity with $k$. Now, the fact that $D_k$ is included in $\bar N_k$ is ensured when $\partial \bar N_k$ is mean convex and one can always deform the Riemannian metric in a small neighbourhood of $\partial \bar N_k$ so that it becomes mean convex. The disk $D_k$ is then included in $\bar N_k$ and is minimal for this new metric  which coincides with the original one in, say, $\bar N_{k-1}\subset N_k$ and this is sufficient for the rest of the argument. Now, let us assume that the sequence $D_k$ converges to a complete minimal surface $\Sigma\subset X$ which, in that situation, is stable. By the result of Schoen and Yau (see \cite{SchoenYau:stable}) this surface is homeomorphic to $\textbf{R}^2$ and the previous argument shows that the number of connected components of $\Sigma\cap N_2$ intersecting $N_1$ is infinite.
By a result of Meeks and Yau (see \cite{Mee-Yau}), each of these components contains a definite amount of area. The contradiction comes from a nice extrinsic version of Cohn-Vossen's Inequality, proved by Jian Wang,
$$\int_\Sigma \kappa (x)dv(x) \leq 2\pi\,,$$
where $\kappa (x)$ is the scalar curvature of $X$ at $x$ and $dv$ is the volume element of the induced metric on $\Sigma$. The original Cohn-Vossen's Inequality is the same, for a complete surface with positive curvature, where $\kappa$ is replaced by the intrinsic curvature. Since on $N_2$ the scalar curvature of $X$ is bounded away from zero, by assumption, this inequality is in contradiction with the infinite area contained in $N_2$. Now, this is by far too naive; indeed, in general $D_k$ does not converge to a complete stable $\Sigma$ but, according to Colding and Minicozzi (see \cite{Col-Min}), to a lamination with complete stable minimal leaves. A variation of the above argument, much more involved, leads to the same contradiction with the extrinsic Cohn-Vossen's Inequality. As mentioned before, in the case when the scalar curvature is non negative, a trick due to J.~Kazdan (see \cite{Kaz}) allows to deform it into a metric with positive scalar curvature. Let us point out that this beautiful and efficient version of Cohn-Vossen's Inequality relies on a smart use of the Jacobi operator for a stable minimal surface (see Theorem 5.10 in \cite{Wang1}).

In the general case, when the open contractible manifold is not genus one, we get again a lamination whose leaves are complete stable minimal surfaces $\Sigma$ homeomorphic to $\textbf{R}^2$ which we call \emph{stable planes}. The properties of these  planes and the way they intersect the $N_k$'s are central in the proofs. The proof of Theorem \ref{theo:Wang} rely on using these stable planes, chosen with specific properties and which are separating $X$ in two connected components, to do some sort of \emph{plane surgery}. Jian wang shows the existence of an uttermost plane $L$ so that, at least, one side of $X\setminus L$, say $X'$, has trivial fundamental group. We could then apply the proof described in \cite{Wang2} except that $X'$ has a boundary. This causes lots of difficulties, for example the stable minimal laminations that one obtains by taking limits of meridian disks may intersect $L$. One has to show that, somehow, $L$ acts as a barrier. It is not the point of this text to go further into the technical issues but we insist on the fact that it is a \emph{tour de force}.

\section{Some Questions}

This short account of the Riemannian geometry of some open and contractible $3$-manifolds is focused on the positive or non negative scalar curvature. Despite these restrictions it is quite difficult to prove results and the techniques involved are definitely sophisticated. 

We did not address any issue concerning a metric whose curvature is bounded above. It is clear, thanks to Cartan-Hadamard Theorem, that any contractible $3$-manifold cannot carry a complete metric with non positive sectional curvature unless it is diffeomorphic to  $\mathbf{R}^{3}$. It turns out that this statement is true even if we relax the regularity. More precisely, D.~Rolfsen proved in \cite{Rolf} that a complete open $\CAT(0)$ manifold of dimension $3$ is homeomorphic to $\mathbf{R}^{3}$. Hence, Whitehead's type manifolds cannot carry any complete $\CAT(0)$ metric. Notice that the above result by  D.~Rolfsen is not any more true in higher dimension (see \cite{D-J}).

In recent years we have witnessed a huge activity around synthetic versions of the notion of Ricci curvature bounded below. New families of metric spaces have been described, such as $\CD (\kappa , N)$ spaces (see \cite{LV}), and lots of results show that they share plenty of nice properties with manifolds whose (standard) Ricci curvature is bounded below by $\kappa$ and whose dimension is (bounded above by) $N$. These notions might not see the details of the local geometry, they rather focus on the geometry in the large; however, let us recall that the assumption in Theorem \ref{cwy} is only at infinity, more precisely outside a compact set. Then, in view of Liu's result (see \cite{Liu}), we are led to ask the following question.

\begin{question}
Does the Whitehead manifold (or any contractible open $3$-manifold not homeomorphic to $\mathbf R^3$) carry a geodesically complete $\CD (0 , 3)$ metric ?
\end{question}

One difficulty is that most of these spaces do not have any  quotient which is a manifold or even an orbifold. We then loose all the tools that group actions could bring into play.
 
Now, going to dimension $4$, there is a family of open spaces which plays a role comparable to the contractible open $3$-manifolds, there are the differentiable structures on $\mathbf R^4$. To our knowledge very little is known about their possible Riemannian geometries and we could dream of a result that would characterise the standard differentiable structure on $\mathbf R^4$ by some of its geometric properties. Some of these exotic $\mathbf R^4$'s, the so-called {\sl large exotic} $\mathbf R^4$, are related to special knots in $\mathbf S^3$ (those which are {\sl topologically slice} but not {\sl smoothly slice}, see the discussion \href{https://mathoverflow.net/questions/42624/slice-knots-and-exotic-mathbb-r4}{here}). It seems possible and plausible that minimal surfaces or hypersurfaces theory will play a role in this $4$-dimensional context too. 

What about a synthetic version of the non negative scalar curvature ?  Some recent works by M.~Gromov seem to pave the way towards such a notion; the interested reader is referred to \cite{Misha1, Misha2} and the preprints posted \href{https://www.ihes.fr/~gromov/category/positivescalarcurvature/}{here}. Other interesting articles along these lines are \cite{ChaoLi1, ChaoLi2}, by Chao Li. This is a very challenging and exciting question.

\section{The Poincaré Conjecture}

The starting point of our interest is the \lq\lq proof " of the Poincaré conjecture given in 1934 by J.~H.~C.~Whitehead in the article \cite{White1}. Roughly speaking the scheme is the following (see \href{https://www.math.unl.edu/~mbrittenham2/ldt/poincare.html}{here}):
\begin{itemize}
   \item let $X$ be  a simply connected closed $3$-manifold then, for any point $p\in X$, $X\setminus\{p\}$ is contractible (by simple considerations left to the reader).
   \item The only contractible $3$-manifold is $\textbf{R}^3$.
   \item The one-point compactification of $\textbf{R}^3$ is $\textbf{S}^3$.
\end{itemize}
In 1935, he realised that the second step was wrong (see \cite{White2}) and constructed in \cite{White3} the first contractible $3$-manifold which is not homeomorphic to $\textbf{R}^3$, the \textsl{Whitehead manifold} which we have denoted by $Wh$ in the previous sections. 

Let us recall that Henri Poincaré never stated this problem as a conjecture but as a question and this distinction is important since, from 1904 until the 70's, there were as many (wrong) proofs of the \lq\lq conjecture" as (wrong) counter-examples provided. Another similar question, but related to open manifolds, could have been asked as follows: is an open manifold which has the same homotopy as $\textbf{R}^3$, that is which is contractible, homeomorphic to $\textbf{R}^3$ ? The Whitehead manifold is a counter-example to this statement and the answer to this question is negative. Nevertheless, the discovery of $3$-manifolds which are contractible but not homeomorphic to $\textbf{R}^3$, such as the Whitehead manifold, opened a wide playground for topologists and geometers. 

The manifold $Wh$ is a simple manifold, namely an open subset of $\textbf{S}^3$ (hence of $\textbf{R}^3$) whose complement is a closed set, called the Whitehead continuum, which looks locally like a product of an interval by a Cantor set. Now, such manifolds, not homeomorphic to $\textbf{R}^3$ but nevertheless  contractible, turn out to be plentiful, as was shown by McMillan (see \cite{McMil}). In fact there exist uncountably many such manifolds whereas there are only a countable family of closed topological manifolds. In \cite{McMil2} (Theorem 1, p. 511), McMillan showed that all contractible open $3$-manifolds satisfying an extra topological assumption are an increasing union of handle bodies, generalising the genus one construction. The assumption that had to be satisfied is that each compact subset of the 3-manifolds considered can be embedded in the 3-sphere $\textbf{S}^3$. Such a
3-manifold was called a W-space. This assumption was dictated by the fact that the Poincaré conjecture was not settled and is now useless. Notice also that $X\setminus\{p\}$ has a topological boundary at infinity which is homeomorphic to $\textbf{S}^2$ (the boundary of a ball centred at $p$), which is clearly not the case for, at least, the Whitehead manifold. Many other nice properties can be proved and the reader is referred to the literature. Finally, it would be nice to have a \lq\lq proof\," along the lines suggested by Whitehead but, at the moment, an idea is still missing.

A general picture of the world of closed $3$-manifolds was given by W.~Thurston in the 70's and is known as Thurston's geometrization conjecture (see \cite{Thur} and the post-Perelman literature). It is a conjectural description of all closed $3$-manifolds and the way to build or decompose them nicely, which yields the Poincaré conjecture as a corollary. The picture is so clear and so nice that the attempts to publish counter-examples to the Poincaré conjecture (almost) stopped. In any case it took then less than 30 years to get G.~Perelman's proof of Thurston's Conjecture (see \cite{Per1, Per2, Per3}) developing a technique introduced by R.~Hamilton in \cite{Ham}.

The beautiful contribution by Thurston is to decompose the closed $3$-manifolds into pieces carrying each a very specific Riemannian geometry, like the case of closed surfaces which are either spherical, flat or hyperbolic but with more (in fact 8) possibilities. Now, what could be a Geometrization conjecture for open $3$-manifolds ? There is little hope to be able to state a reasonable conjecture. For example, the starting point of Thurston's Geometrization Conjecture are the Kneser-Milnor and the Jaco-Shalen and Johannson decompositions of closed $3$-manifolds none of which is true for open manifolds (see \cite{S-T} and \cite{Mai}). Therefore, the only sensible thing to do is to work out examples and, in this spirit, the purpose of the previous sections is to describe an attempt to understand the geometry of the family of contractible open $3$-manifolds.

\bibliographystyle{alpha}
\bibliography{bibliography-Chern}
\end{document}